\documentclass[12 pt, a4 paper]{article}
\usepackage{amssymb,amsmath,graphicx,enumerate}
\usepackage{longtable,mathrsfs}

\numberwithin{equation}{section}
\newtheorem{theorem}{Theorem}[section]

\newtheorem{definition}{Definition}[section]

\newtheorem{prop}{Proposition}[section]
\allowdisplaybreaks

\begin{document}
\begin{center}\begin{large}
P$\acute{o}$lya-Szeg$\ddot{o}$ fractional inequality via Hadamard fractional integral
 \end{large}\end{center}
\begin{center}
                 $Vaijanath  \, L. Chinchane ^{1},\,\,\,   Deepak \, B. Pachpatte ^{2},\,\,\,   Asha \, B. Nale^{3}$

$^{1}$Department of Mathematics,\\
Deogiri Institute of Engineering and Management\\
Studies Aurangabad-431005, INDIA\\
E-mail ID: chinchane85@gmail.com\\
$^{2}$Department of Mathematics,\\
Dr. Babasaheb Ambedkar Marathwada University,\\
 Aurangabad-431 004, INDIA.\\
E-mail ID: pachpatte@gmail.com\\
$^{3}$Department of Mathematics,\\
Dr. Babasaheb Ambedkar Marathwada University,\\
 Aurangabad-431 004, INDIA.\\
E-mail ID: ashabnale@gmail.com

\end{center}
\begin{abstract}
 In the present paper, we establish some new fractional integral inequalities similar to P$\acute{o}$lya-Szeg$\ddot{o}$ integral inequality and fractional inequality related to Minkowsky inequality by using the Hadamard fractional integral operator. Also, we discuss few spacial cases of these inequalities.
\end{abstract}
\textbf{Keywords:} P$\acute{o}$lya-Szeg$\ddot{o}$ inequality, Hadamard fractional integral and integral inequality.\\
\textbf{Mathematics Subject Classification :} 26D10, 26A33.\\
\section{Introduction}
\paragraph{}It is well known truth that inequalities have a big important in development of many branches of mathematics and other fields of sciences.
In 1935, G. Gr$\ddot{u}$ss proved the following classical integral inequality, see \cite{GG1}
\begin{equation}
\begin{split}
&\left|\frac{1}{b-a}\int^{b}_{a}f(x)g(x)dx-\left(\frac{1}{b-a}\int^{b}_{a}f(x)dx\right)\left(\frac{1}{b-a}\int^{b}_{a}g(x)dx\right)\right|\\
&\leq \frac{(M-m)(P-p)}{4},
\end{split}
\end{equation}
provided that $f$ and $g$ are two integrable functions on $[a, b]$ and satisfy the condition
\begin{equation}
m\leq f(x)\leq M,\,\,\, p\leq g(x)\leq P;\,\,\,m, M, p, P \in R, \, x\in[a, b].
\end{equation}
In 1925, P$\acute{o}$lya-Szeg$\ddot{o}$ proved the following inequality, see [13]
\begin{equation}
\frac{\int_{a}^{b}f^{2}(x)dx\int_{a}^{b}g^{2}(x)dx}{\left(\int_{a}^{b}f(x)dx\int_{a}^{b}g(x)dx\right)^{2}}\leq \frac{1}{4} \left( \sqrt\frac{MN}{mn}+\sqrt{\frac{mn}{MN}}\right)^{2},
\end{equation}
and in [9] Dargomir and Diamond proved the following inequality
\begin{equation}
\begin{split}
&\left|\frac{1}{b-a}\int^{b}_{a}f(x)g(x)dx-\left(\frac{1}{b-a}\int^{b}_{a}f(x)dx\right)\left(\frac{1}{b-a}\int^{b}_{a}g(x)dx\right)\right|\\
&\leq \frac{(M-m)(P-p)}{4(b-a)^{2}\sqrt{mMpP}}\int^{b}_{a}f(x)dx\int^{b}_{a}g(x)dx,
\end{split}
\end{equation}
provided that $f$ and $g$ are two integrable functions on $[a, b]$ and satisfy the condition
\begin{equation}
0<m\leq f(x)\leq M <\infty,\,\,\, 0<p\leq g(x)\leq P<\infty;\,\,\,m, M, p, P \in R, \, x\in[a, b].
\end{equation}
During last two decades many researchers have work on the fractional integral inequalities by using Riemann-Liouville and Saigo integral operator, see [7, 8, 10, 14, 15].
In [14] Shilpi Jain and et al. have proved the P$\acute{o}$lya-Szeg$\ddot{o}$ inequality using Riemann-Liouville fractional integral operator which is as follows
\begin{theorem} let $f$ and $g$ be two integrable functions on $[0, \infty]$. Assume that there exist four positive integrable functions $u_{1}$,$u_{2}$, $v_{1}$ and $v_{2}$ on $[0, \infty]$ such that
\begin{equation*}
0<u_{1}(\tau)\leq f(\tau)\leq u_{2}(\tau), \, \, 0<v_{1}(\tau)\leq g(\tau)\leq v_{2}(\tau).
\end{equation*}
Then for $t>0$ and $\alpha, \beta > 0$ the following inequalities hold:
\begin{equation}
\frac{I_{t}^{\alpha,\beta,\eta,\mu}\{v_{1}v_{2}f^{2}\}(t)I_{t}^{\alpha,\beta,\eta,\mu}\{u_{1}u_{2}g^{2}\}(t)}{I_{t}^{\alpha,\beta,\eta,\mu}\{(v_{1}u_{1}+v_{2}u_{2})fg\}(t)}\leq \frac{1}{4}.
\end{equation}
And
\begin{equation}
\frac{I_{t}^{\alpha,\beta,\eta,\mu}\{u_{1}u_{2}\}(t)I_{t}^{\gamma,\delta,\zeta,\nu}\{v_{1}v_{2}\}(t)I_{t}^{\alpha,\beta,\eta,\mu}\{f^{2}\}(t)I_{t}^{\gamma,\delta,\zeta,\nu}\{g^{2}\}(t)}{\left(I_{t}^{\alpha,\beta,\eta,\mu}\{u_{1}f(t)\}I_{t}^{\gamma,\delta,\zeta,\nu}\{v_{1}g(t)\}I_{t}^{\alpha,\beta,\eta,\mu}\{u_{2}f(t)\}I_{t}^{\gamma,\delta,\zeta,\nu}\{v_{2}g\}(t)\right)^{2}}\leq \frac{1}{4}.
\end{equation}
\end{theorem}
Where $I_{t}^{\alpha,\beta,\eta,\mu}\{.\}$ denote the generalized fractional integral (in terms of the Gauss hypergeometric function)of order $\alpha$ for a real valued continuous function.
\paragraph{} In literature few results are available in which Hadamard fractional integral operator has been used [2-6]. Motivated by the above work this paper we have established some new results for P$\acute{o}$lya- Szeg$\ddot{o}$ inequality and some other inequalities using Hadamard fractional integral.
\section{Preliminaries}
\paragraph{}Here, we present some basic definitions of Hadamard derivative and integral as given in [1, p.159-171].
\begin{definition} \cite{AHJ}
The  Hadamard fractional integral of order $\alpha \in R^{+}$ of function  $f(x)$, for all $x>1$ is defined as,
\begin{equation}
_{H}D_{1,x}^{-\alpha}f(x)= \frac{1}{\Gamma (\alpha)}\int_{1}^{x}\ln(\frac {x}{t})^{\alpha-1} f(t) \frac{dt}{t},
\end{equation}
where $\Gamma (\alpha)= \int_0^\infty e^{-u}u^{\alpha-1} du$.
\end{definition}
\begin{definition}$[1]$
The Hadmard fractional derivative of order $\alpha \in [n-1,n)$, $n\in Z^{+}$, of function $f(x)$ is given as follows
\begin{equation}
_{H}D_{1,x}^{\alpha}f(x)= \frac{1}{\Gamma (n-\alpha)}(x \frac{d}{dx})^{n}\int_{1}^{x}\ln(\frac {x}{t})^{n-\alpha-1} f(t) \frac{dt}{t}.
\end{equation}
\end{definition}
\begin{prop}$[1]$
If $0<\alpha <1$, the following relation hold
\begin{equation}
{_{H}}D_{1,x}^{-\alpha}(\ln{x})^{\beta-1}= \frac{\Gamma (\beta)}{ \Gamma(\beta+\alpha)}(\ln{x})^{\beta+\alpha-1},
\end{equation}
\begin{equation}
_{H}D_{1,x}^{\alpha}(\ln{x})^{\beta-1}= \frac{\Gamma (\beta)}{ \Gamma(\beta-\alpha)}(\ln{x})^{\beta-\alpha-1},
\end{equation}
\end{prop}
respectively.
For the convenience of establishing the result, we give the semigroup property is as follows,
\begin{equation}
(_{H}D_{1,x}^{-\alpha})(_{H}D_{1,x}^{-\beta})f(x)= _{H}D_{1,x}^{-(\alpha+\beta)}f(x).
\end{equation}
Also some details of fractional Hadamard calculus are given in the book A.A.Kilbas et al. [12], and in book of S.G.Samko et al. [16].
\section{ Fractional P$\acute{o}$lya-Szeg$\ddot{o}$ inequality}
\paragraph{} Here, we give our first main theorem.
\begin{theorem} let $x$ and $y$ be two integrable functions on $[1, \infty]$. Assume that there exist four positive integrable functions $u_{1}$,$u_{2}$, $v_{1}$ and $v_{2}$ on $[1, \infty]$ such that
\begin{equation*}
0<u_{1}(\tau)\leq x(\tau)\leq u_{2}(\tau), \, \, 0<v_{1}(\tau)\leq y(\tau)\leq v_{2} (\tau), \, \, (\tau \in (0,t)], t>0).
\end{equation*}
Then for $t>0$ and $\alpha >0$, the following inequality hold
\begin{equation}
\frac{_{H}D_{1,t}^{-\alpha}\{v_{1}v_{2}x^{2}\}(t)_{H}D_{1,t}^{-\alpha}\{u_{1}u_{2}y^{2}\}(t)}{_{H}D_{1,t}^{-\alpha}\{(v_{1}u_{1}+v_{2}u_{2})xy\}(t)}\leq \frac{1}{4}.
\end{equation}
\end{theorem}
\textbf{Proof:-} To prove (3.1), since $\tau \in (0,t)$ and $t>0$ we have
\begin{equation}
\frac{x(\tau)}{y(\tau)}\leq \frac{u_{2}(\tau)}{v_{1}(\tau)},
\end{equation}
which implies that
\begin{equation}
\left(\frac{u_{2}(\tau)}{v_{1}(\tau)}-\frac{x(\tau)}{y(\tau)}\right)\geq 0,
\end{equation}
also we have,
\begin{equation}
\left(\frac{x(\tau)}{y(\tau)}-\frac{u_{1}(\tau)}{v_{2}(\tau)}\right)\geq 0,
\end{equation}
Multiplying equation (3.3) and (3.4), we have
\begin{equation}
\left(\frac{u_{2}(\tau)}{v_{1}(\tau)}-\frac{x(\tau)}{y(\tau)}\right)\left(\frac{x(\tau)}{y(\tau)}-\frac{u_{1}(\tau)}{v_{2}(\tau)}\right)\geq 0,
\end{equation}
\begin{equation}
\left(\frac{u_{2}(\tau)}{v_{1}(\tau)}-\frac{x(\tau)}{y(\tau)}\right)\frac{x(\tau)}{y(\tau)}-\left(\frac{u_{2}(\tau)}{v_{1}(\tau)}-\frac{x(\tau)}{y(\tau)}\right)\frac{u_{1}(\tau)}{v_{2}(\tau)}\geq 0,
\end{equation}
it follows that
\begin{equation}
\left(\frac{u_{2}(\tau)}{v_{1}(\tau)}+\frac{u_{1}(\tau)}{v_{2}(\tau)}\right)\frac{x(\tau)}{y(\tau)}\geq \frac{x^{2}(\tau)}{y^{2}(\tau)}+\frac{u_{1}(\tau)u_{2}(\tau)}{v_{1}(\tau)v_{2}(\tau)},
\end{equation}
after some calculation we obtain,
\begin{equation}
\left[u_{1}(\tau)v_{1}(\tau)+u_{2}(\tau)v_{2}(\tau)\right]x(\tau)y(\tau)\geq v_{1}(\tau)v_{2}(\tau)x^{2}(\tau)+u_{2}(\tau)v_{2}(\tau)y^{2}(\tau),
\end{equation}
multiplying equation by (3.8) by $ \frac{({\ln(\frac {t}{\tau})})^{\alpha-1}}{\tau \Gamma(\alpha)}$, we obtain
\begin{equation}
\begin{split}
&\frac{({\ln(\frac {t}{\tau})})^{\alpha-1}}{\tau \Gamma(\alpha)}\left[u_{1}(\tau)v_{1}(\tau)+u_{2}(\tau)v_{2}(\tau)\right]x(\tau)y(\tau)\\
&\geq \frac{({\ln(\frac {t}{\tau})})^{\alpha-1}}{\tau \Gamma(\alpha)}v_{1}(\tau)v_{2}(\tau)x^{2}(\tau)+\frac{({\ln(\frac {t}{\tau})})^{\alpha-1}}{\tau \Gamma(\alpha)}u_{2}(\tau)v_{2}(\tau)y^{2}(\tau).
\end{split}
\end{equation}
Integrate the equation (3.9) on both side with respective $\tau$ from $1$ to $t$, we get
\begin{equation}
_{H}D_{1,t}^{-\alpha}\left[u_{1}v_{1}+u_{2}v_{2}\right]xy(t)\geq _{H}D_{1,t}^{-\alpha}v_{1}v_{2}x^{2}(t)+_{H}D_{1,t}^{-\alpha}u_{2}v_{2}y^{2}(t).
\end{equation}
By using elementary inequalities $a+b\geq \sqrt{ab},$\,  where $a, b \in R^{+} $, we have
\begin{equation}
_{H}D_{1,t}^{-\alpha}\left[u_{1}v_{1}+u_{2}v_{2}\right]xy(t)\geq 2\sqrt{_{H}D_{1,t}^{-\alpha}v_{1}v_{2}x^{2}(t)\,_{H}D_{1,t}^{-\alpha}u_{2}v_{2}y^{2}(t)},
\end{equation}
squaring both side of (3.11)
\begin{equation}
\left(_{H}D_{1,t}^{-\alpha}\left[u_{1}v_{1}+u_{2}v_{2}\right]xy(t)\right)^{2}\geq 4\left(_{H}D_{1,t}^{-\alpha}v_{1}v_{2}x^{2}(t)\,_{H}D_{1,t}^{-\alpha}u_{2}v_{2}y^{2}(t)\right).
\end{equation}
It follows that
\begin{equation}
_{H}D_{1,t}^{-\alpha}v_{1}v_{2}x^{2}(t)\,_{H}D_{1,t}^{-\alpha}u_{2}v_{2}y^{2}(t) \leq \frac{1}{4}\left(_{H}D_{1,t}^{-\alpha}\left[u_{1}v_{1}+u_{2}v_{2}\right]xy(t)\right)^{2},
\end{equation}
which gives the required inequality (3.1).
\begin{theorem} let $x$ and $y$ be two integrable functions on $[1, \infty]$. Assume that there exist four positive integrable functions $u_{1}$,$u_{2}$, $v_{1}$ and $v_{2}$ on $[1, \infty]$ such that
\begin{equation*}
0<u_{1}(\tau)\leq x(\tau)\leq u_{2}(\tau), \, \, 0<v_{1}(\rho)\leq y(\rho)\leq v_{2} (\rho), \, \, (\tau, \rho \in [0,t], t>0).
\end{equation*}
Then for $t>0$ and $\alpha >0, \beta>0$, the following inequality hold
\begin{equation}
\frac{_{H}D_{1,t}^{-\alpha}\{u_{1}u_{2}\}(t)_{H}D_{1,t}^{-\beta}\{v_{1}v_{2}\}(t)_{H}D_{1,t}^{-\alpha}\{x^{2}\}(t)_{H}D_{1,t}^{-\beta}\{y^{2}\}(t)}{\left(_{H}D_{1,t}^{-\alpha}\{(u_{1}x)\}(t)_{H}D_{1,t}^{-\beta}\{(v_{1}y)\}(t)+_{H}D_{1,t}^{-\alpha}\{(u_{2}x)\}(t)_{H}D_{1,t}^{-\beta}\{(v_{2}y)\}(t)\right)^{2}}\leq \frac{1}{4}.
\end{equation}
\end{theorem}
\textbf{Proof:-} To prove (3.14), since $\tau, \rho \in (0,t]$ and $t>0$ we have
\begin{equation}
\frac{x(\tau)}{y(\rho)}\leq \frac{u_{2}(\tau)}{v_{1}(\rho)},
\end{equation}
which implies that
\begin{equation}
\left(\frac{u_{2}(\tau)}{v_{1}(\rho)}-\frac{x(\tau)}{y(\rho)}\right)\geq 0,
\end{equation}
also we have,
\begin{equation}
\left(\frac{x(\tau)}{y(\rho)}-\frac{u_{1}(\tau)}{v_{2}(\rho)}\right)\geq 0,
\end{equation}
Multiplying equation (3.16) and (3.17), we have
\begin{equation}
\left(\frac{u_{2}(\tau)}{v_{1}(\rho)}-\frac{x(\tau)}{y(\rho)}\right)\left(\frac{x(\tau)}{y(\rho)}-\frac{u_{1}(\tau)}{v_{2}(\rho)}\right)\geq 0,
\end{equation}
\begin{equation}
\left(\frac{u_{2}(\tau)}{v_{1}(\rho)}-\frac{x(\tau)}{y(\rho)}\right)\frac{x(\tau)}{y(\rho)}-\left(\frac{u_{2}(\tau)}{v_{1}(\rho)}-\frac{x(\tau)}{y(\rho)}\right)\frac{u_{1}(\tau)}{v_{2}(\rho)}\geq 0,
\end{equation}
it follows that
\begin{equation}
\left(\frac{u_{2}(\tau)}{v_{1}(\rho)}+\frac{u_{1}(\tau)}{v_{2}(\rho)}\right)\frac{x(\tau)}{y(\rho)}\geq \frac{x^{2}(\tau)}{y^{2}(\rho)}+\frac{u_{1}(\tau)u_{2}(\tau)}{v_{1}(\rho)v_{2}(\rho)}.
\end{equation}
Multiplying both side of equation (3.20) by $v_{1}(\rho)v_{}(\rho)y^{2}(\rho)$, we obtain
\begin{equation}
u_{1}(\tau)x(\tau)v_{1}(\rho)g(\rho)+u_{2}(\tau)x(\tau)v_{1}(\rho)g(\rho) \geq v_{1}(\rho)v_{2}(\rho)x^{2}(\tau)+u_{1}(\tau)u_{2}(\tau)y^{2}(\rho),
\end{equation}
multiplying both side of (3.21) by $ \frac{({\ln(\frac {t}{\tau})})^{\alpha-1}}{\tau \Gamma(\alpha)}$, which is positive because $\tau \in(0,t)$, $t>0$, then integrate the resulting identity with respect to $\tau$ from $1$ to $t$, we get
\begin{equation}
\begin{split}
&v_{1}(\rho)g(\rho)_{H}D_{1,t}^{-\alpha}u_{1}x(t)+v_{1}(\rho)g(\rho)_{H}D_{1,t}^{-\alpha}u_{2}x(t)\\
 &\geq v_{1}(\rho)v_{2}(\rho)_{H}D_{1,t}^{-\alpha}x^{2}(t)+y^{2}(\rho)_{H}D_{1,t}^{-\alpha}u_{1}u_{2}(t),
\end{split}
\end{equation}
now, again multiplying both side of (3.22) by $ \frac{({\ln(\frac {t}{\rho})})^{\alpha-1}}{ \rho \Gamma(\alpha)}$, which is positive because $\rho\in(0,t)$, $t>0$, then integrate the resulting identity with respect to $\rho$ from $1$ to $t$, we get
\begin{equation}
\begin{split}
&_{H}D_{1,t}^{-\beta}v_{1}g(t)_{H}D_{1,t}^{-\alpha}u_{1}x(t)+_{H}D_{1,t}^{-\beta}v_{1}g(t)_{H}D_{1,t}^{-\alpha}u_{2}x(t)\\
 &\geq _{H}D_{1,t}^{-\beta}v_{1}v_{2}(t)_{H}D_{1,t}^{-\alpha}x^{2}(t)+_{H}D_{1,t}^{-\beta}y^{2}(t)_{H}D_{1,t}^{-\alpha}u_{1}u_{2}(t).
\end{split}
\end{equation}
By using elementary inequalities $a+b\geq \sqrt{ab},$\,  where $a, b \in R^{+} $, we have
\begin{equation}
\begin{split}
&_{H}D_{1,t}^{-\beta}v_{1}g(t)_{H}D_{1,t}^{-\alpha}u_{1}x(t)+_{H}D_{1,t}^{-\beta}v_{1}g(t)_{H}D_{1,t}^{-\alpha}u_{2}x(t)\\
&\geq 2\sqrt{_{H}D_{1,t}^{-\beta}v_{1}v_{2}(t)_{H}D_{1,t}^{-\alpha}x^{2}(t)_{H}D_{1,t}^{-\beta}y^{2}(t)_{H}D_{1,t}^{-\alpha}u_{1}u_{2}(t)},
\end{split}
\end{equation}
which give the required inequality (3.14). This complete the proof.
\begin{theorem} let $x$ and $y$ be two integrable functions on $[1, \infty]$. Assume that there exist four positive integrable functions $u_{1}$,$u_{2}$, $v_{1}$ and $v_{2}$ on $[1, \infty]$ such that
\begin{equation*}
0<u_{1}(\tau)\leq x(\tau)\leq u_{2}(\tau), \, \,  0<v_{1}(\rho)\leq y(\rho)\leq v_{2} (\rho), \, \, (\tau, \rho \in [0,t], t>0).
\end{equation*}
Then for $t>0$ and $\alpha, \beta>0$, the following inequality hold
\begin{equation}
_{H}D_{1,t}^{-\alpha}\{x^{2}\}(t)_{H}D_{1,t}^{-\beta}\{y^{2}\}(t)\leq _{H}D_{1,t}^{-\alpha}\{\frac{u_{2}xy}{v_{1}}\}(t)_{H}D_{1,t}^{-\beta}\{\frac{v_{2}xy}{u_{1}}\}(t).
\end{equation}
\end{theorem}
\textbf{Proof:-} Multiplying (3.2) by $x(\tau)$, we obtain
\begin{equation}
 x^{2}(\tau)\leq \frac{u_{2}(\tau)}{v_{1}(\tau)}x(\tau)y(\tau),
\end{equation}
Multiplying both side of (3.26) by $ \frac{({\ln(\frac {t}{\tau})})^{\alpha-1}}{\tau \Gamma(\alpha)}$, which is positive because $\tau \in(0,t)$, $t>0$, then integrate the resulting identity with respect to $\tau$ from $1$ to $t$, we get
\begin{equation}
 _{H}D_{1,t}^{-\alpha}x^{2}(t)\leq _{H}D_{1,t}^{-\alpha}\{\frac{u_{2}xy}{v_{1}}\}(t).
\end{equation}
Analogously, we obtain
\begin{equation}
 _{H}D_{1,t}^{-\alpha}y^{2}(t)\leq _{H}D_{1,t}^{-\alpha}\{\frac{v_{2}xy}{u_{1}}\}(t),
\end{equation}
multiplying the inequality (3.27) and (3.28), we obtain the required inequality (3.25). This complete the proof.
\paragraph{} Here we present some special case of above theorem which is as below
\begin{prop} let $x$ and $y$ be two integrable functions on $[1, \infty]$ such that
\begin{equation}
0<m\leq x(\tau)\leq M<\infty, \, \, 0<n \leq y(\tau)\leq N <\infty, \, \, (\tau \in [0,t], t>0).
\end{equation}
Then for $t>0$ and $\alpha >0$, the following inequality hold
\begin{equation}
\frac{_{H}D_{1,t}^{-\alpha}\{x^{2}\}(t)_{H}D_{1,t}^{-\alpha}\{y^{2}\}(t)}{(_{H}D_{1,t}^{-\alpha}\{xy\}(t))^{2}}\leq \frac{1}{4}\left(\sqrt{\frac{mn}{MN}}+\sqrt{\frac{MN}{mn}}\right)^{2}.
\end{equation}
\end{prop}
\begin{prop} let $x$ and $y$ be two integrable functions on $[1, \infty]$ satisfies condition .
\begin{equation}
0<m\leq x(\tau)\leq M<\infty, \, \, 0<n \leq y(\rho)\leq N <\infty, \, \, (\tau, \rho \in [0,t], t>0).
\end{equation}
Then for $t>0$ and $\alpha, \beta >0$, the following inequality hold
\begin{equation}
\frac{(\ln t)^{\alpha+\beta}}{\Gamma(\alpha+1)\Gamma(\beta+1)}\frac{_{H}D_{1,t}^{-\alpha}\{x^{2}\}(t)_{H}D_{1,t}^{-\beta}\{y^{2}\}(t)}{(_{H}D_{1,t}^{-\alpha}\{x\}(t)(_{H}D_{1,t}^{-\beta}\{y\}(t))^{2}}\leq \frac{1}{4}\left(\sqrt{\frac{mn}{MN}}+\sqrt{\frac{MN}{mn}}\right)^{2}.
\end{equation}
\end{prop}
\begin{prop} let $x$ and $y$ be two integrable functions on $[1, \infty]$ satisfies condition (3.29).
Then for $t>0$ and $\alpha, \beta >0$, we have
\begin{equation}
\frac{_{H}D_{1,t}^{-\alpha}\{x^{2}\}(t)_{H}D_{1,t}^{-\beta}\{y^{2}\}(t)}{(_{H}D_{1,t}^{-\alpha}\{xy\}(t)(_{H}D_{1,t}^{-\beta}\{xy\}(t))^{2}}\leq \frac{MN}{mn}.
\end{equation}
\end{prop}
Here, we present the fractional integral inequality related to Minkowsky inequality as follows
\begin{theorem} let $x$ and $y$ be two integrable functions on $[1, \infty]$ such that $\frac{1}{p}+\frac{1}{q}=1, p>1,$ and $0<m<\frac{x(\tau)}{y(\tau)}<M, \tau \in (0,t), t>0.$ Then for all $\alpha>0,$ we have
\begin{equation}
_{H}D_{1,t}^{-\alpha}\{xy\}(t)\leq \frac{2^{p-1}M^{p}}{p(M+1)^{p}}\left(_{H}D_{1,t}^{-\alpha}[x^{p}+y^{p}](t)\right)+\frac{2^{q-1}}{q(m+1)^{q}}\left(_{H}D_{1,t}^{-\alpha}[x^{q}+y^{q}](t)\right).
\end{equation}
\end{theorem}
\textbf{Proof:-} Since,  $\frac{x(\tau)}{y(\tau)}<M, \tau \in (0,t), t>0,$ we have
\begin{equation}
(M+1)x(\tau)\leq M(x+y)(\tau),
\end{equation}
Taking p th power on both side and multiplying resulting identity by $ \frac{({\ln(\frac {t}{\tau})})^{\alpha-1}}{\tau \Gamma(\alpha)}$, we obtain
\begin{equation}
\frac{({\ln(\frac {t}{\tau})})^{\alpha-1}}{\tau \Gamma(\alpha)}(M+1)^{p}x^{p}(\tau)\leq  \frac{({\ln(\frac {t}{\tau})})^{\alpha-1}}{\tau \Gamma(\alpha)}M^{p}(x+y)^{p}(\tau),
\end{equation}
integrate the equation (3.36) on both side with respective $\tau$ from $1$ to $t$, we get
\begin{equation}
\frac{1}{\Gamma(\alpha)} \int^{t}_{1}({\ln(\frac {t}{\tau})})^{\alpha-1}(M+1)^{p}x^{p}(\tau)\frac{d\tau}{\tau} \leq  \frac{1}{\Gamma(\alpha)} \int^{t}_{1}({\ln(\frac {t}{\tau})})^{\alpha-1}M^{p}(x+y)^{p}(\tau)\frac{d\tau}{\tau},
\end{equation}
therefore,
\begin{equation}
_{H}D_{1,t}^{-\alpha}[x^{p}(t)]\leq \frac{M^{p}}{(M+1)^{p}}_{H}D_{1,t}^{-\alpha}[(x+y)^{p}(t)],
\end{equation}
on other hand, $0<m<\frac{x(\tau)}{y(\tau)}, \tau \in (0,t), t>0, we can write$
\begin{equation}
(m+1)y(\tau)\leq (x+y)(\tau),
\end{equation}
therefore,
\begin{equation}
\frac{1}{\Gamma(\alpha)} \int^{t}_{1}({\ln(\frac {t}{\tau})})^{\alpha-1}(m+1)^{q}y^{q}(\tau)\frac{d\tau}{\tau} \leq  \frac{1}{\Gamma(\alpha)} \int^{t}_{1}({\ln(\frac {t}{\tau})})^{\alpha-1}(x+y)^{q}(\tau)\frac{d\tau}{\tau},
\end{equation}
consequently, we have
\begin{equation}
_{H}D_{1,t}^{-\alpha}[y^{q}(t)]\leq \frac{1}{(m+1)^{q}}_{H}D_{1,t}^{-\alpha}[(x+y)^{q}(t)].
\end{equation}
Now, by using Young inequality, we have
\begin{equation}
[x(\tau)y(\tau)]\leq \frac{x^{p}(\tau)}{p}+\frac{y^{q}(\tau)}{q}.
\end{equation}
Multiplying both side of (3.42) by $ \frac{({\ln(\frac {t}{\tau})})^{\alpha-1}}{\tau \Gamma(\alpha)}$, which is positive because $\tau \in(0,t)$, $t>0$, then integrate the resulting identity with respect to $\tau$ from $1$ to $t$, we get
\begin{equation}
_{H}D_{1,t}^{-\alpha}[x(t)y(t))]\leq \frac{1}{p}\,_{H}D_{1,t}^{-\alpha}[x^{p}(t)]+\frac{1}{q}\,_{H}D_{1,t}^{-\alpha}[y^{q}(t)],
\end{equation}
from equation (3.38), (3.41) and (3.43) we obtain
\begin{equation}
_{H}D_{1,t}^{-\alpha}[x(t)y(t))]\leq \frac{M^{p}}{p(M+1)^{p}}\,_{H}D_{1,t}^{-\alpha}[(x+y)^{p}(t)]+\frac{1}{q(m+1)^{q}}\,_{H}D_{1,t}^{-\alpha}[(x+y)^{q}(t)],
\end{equation}
now using the inequality $(a+b)^{r}\leq 2^{r-1}(a^{r}+b^{r}), r>1, a,b \geq 0,$ we have
\begin{equation}
_{H}D_{1,t}^{-\alpha}[(x+y)^{p}(t)] \leq 2^{p-1}\,_{H}D_{1,t}^{-\alpha}[(x^{p}+y^{p})(t)],
\end{equation}
and
\begin{equation}
_{H}D_{1,t}^{-\alpha}[(x+y)^{q}(t)] \leq 2^{q-1}\,_{H}D_{1,t}^{-\alpha}[(x^{q}+y^{q})(t)].
\end{equation}
Injecting (3.45), (3.46) in (3.44) we get required inequality (3.34). this complete the proof.\\
\textbf{Concluding remark}\\
It is conclude that the results give some contributions to the theory of integral inequalities and fractional calculus. Moreover, they are expected to lead to some applications for establishing uniqueness of solutions in fractional  differential equations.\\
\textbf{Competing interests}\\
The authors declare that they have no competing interests.\\
\textbf{Authors’ contributions}\\
All authors have equal contributions. All authors read and approved the final manuscript.
		
\end{document}